\documentclass[11pt,a4paper]{article}
\usepackage{cite}
\usepackage{amsmath, amsthm, amssymb,slashed}
\usepackage{ifpdf}
\ifpdf
  \usepackage[pdftex]{graphicx}
  \usepackage{epstopdf}
\else
  \usepackage[dvips]{graphicx}
\fi
\parskip=\baselineskip
\def\Bbb{\mathbb}
\def\Tr{{\rm Tr}}
\def\16{{\bf 16}}
\def\1{{\bf 1}}
\def\2{{\bf 2}}
\def\3{{\bf 3}}
\def\4{{\bf 4}}
 \def\Sp{{\mathrm{Sp}}}
 \def\Spin{{\mathrm{Spin}}}
 \def\SU{{\mathrm{SU}}}
 \def\SO{{\mathrm{SO}}}
\def\bar{\overline}

\def\R{{\Bbb{R}}}\def\Z{{\Bbb{Z}}}
\def\N{{\mathcal N}}

\font\teneurm=eurm10 \font\seveneurm=eurm7 \font\fiveeurm=eurm5
\newfam\eurmfam
\textfont\eurmfam=\teneurm \scriptfont\eurmfam=\seveneurm
\scriptscriptfont\eurmfam=\fiveeurm

\font\teneusm=eusm10 \font\seveneusm=eusm7 \font\fiveeusm=eusm5
\newfam\eusmfam
\textfont\eusmfam=\teneusm \scriptfont\eusmfam=\seveneusm
\scriptscriptfont\eusmfam=\fiveeusm

\font\tencmmib=cmmib10 \skewchar\tencmmib='177
\font\sevencmmib=cmmib7 \skewchar\sevencmmib='177
\font\fivecmmib=cmmib5 \skewchar\fivecmmib='177
\newfam\cmmibfam
\textfont\cmmibfam=\tencmmib \scriptfont\cmmibfam=\sevencmmib
\scriptscriptfont\cmmibfam=\fivecmmib

\numberwithin{equation}{section}

\def\d{\mathrm d}
\def\C{{\Bbb C}}
\def\Z{{\Bbb Z}}
\def\A{{\mathcal A}}
\def\bar{\overline}
\def\W{{\mathcal W}}
\begin{document}
\begin{titlepage}
\begin{flushright}

\end{flushright}
\vskip 1.5in
\begin{center}
{\bf\Large{Two Lectures On The Jones Polynomial}}
\vskip.1cm{\bf\Large{ and Khovanov Homology}}
\vskip
0.5cm {Edward Witten} \vskip 0.05in {\small{ \textit{School of
Natural Sciences, Institute for Advanced Study}\vskip -.4cm
{\textit{Einstein Drive, Princeton, NJ 08540 USA}}}
}
\end{center}
\vskip 0.5in
\baselineskip 16pt
\begin{abstract}
In the first of these two lectures, I describe a gauge theory approach to understanding quantum knot
invariants as Laurent polynomials in a complex variable $q$.  The two main steps are to reinterpret three-dimensional
Chern-Simons gauge theory in four dimensional terms and then to apply electric-magnetic duality.   The variable $q$
is associated to instanton number in the dual description in four dimensions. In the second lecture,
I describe how Khovanov homology can emerge upon adding a fifth dimension.  (Based on lectures presented
at the Clay Research Conference at Oxford University, and also at the Galileo Galilei Institute in Florence, the University of Milan,
Harvard University, and the University of Pennsylvania.)
\end{abstract}
\date{September, 2011}
\end{titlepage}
\def\Hom{\mathrm{Hom}}

\def\U{{\mathcal U}}
\section{Lecture One}

The Jones polynomial is a celebrated invariant of a knot (or link) in ordinary three-dimensional
space, originally discovered by V. F. R. Jones roughly thirty years ago as an offshoot of his work on
von Neumann algebras  \cite{Jones}.  Many descriptions and generalizations of the Jones polynomial were discovered in the years
immediately after Jones's work.  They more or less all involved statistical mechanics or two-dimensional mathematical physics in one way or another -- for example, Jones's 
original work involved Temperley-Lieb algebras of statistical mechanics.  I do not want to assume that the Jones polynomial is familiar to everyone, so I will explain one
of the original definitions.

\begin{figure}[ht]
 \begin{center}
   \includegraphics[width=1.3in]{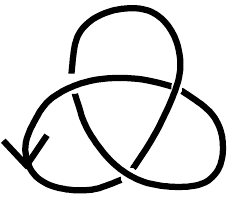}
 \end{center}
\caption{\small A knot in $\R^3$ -- in this case a trefoil knot  -- projected to the plane $\R^2$ in a way that gives an immersion with only simple
crossings and such that the height function (the vertical coordinate in the figure) has only simple local maxima and minima.  In this example, there are three crossings
and two local minima and maxima, making a total of $3+2+2=7$ exceptional points.  Omitting those points divides
the knot into 7 pieces that can be labeled by symbols $+$ or $-$, so the vertex model expresses the Jones
polynomial of the trefoil knot as a sum of $2^7$ terms. \label{projection}}
\end{figure}

\begin{figure}
 \begin{center}
   \includegraphics[width=6in]{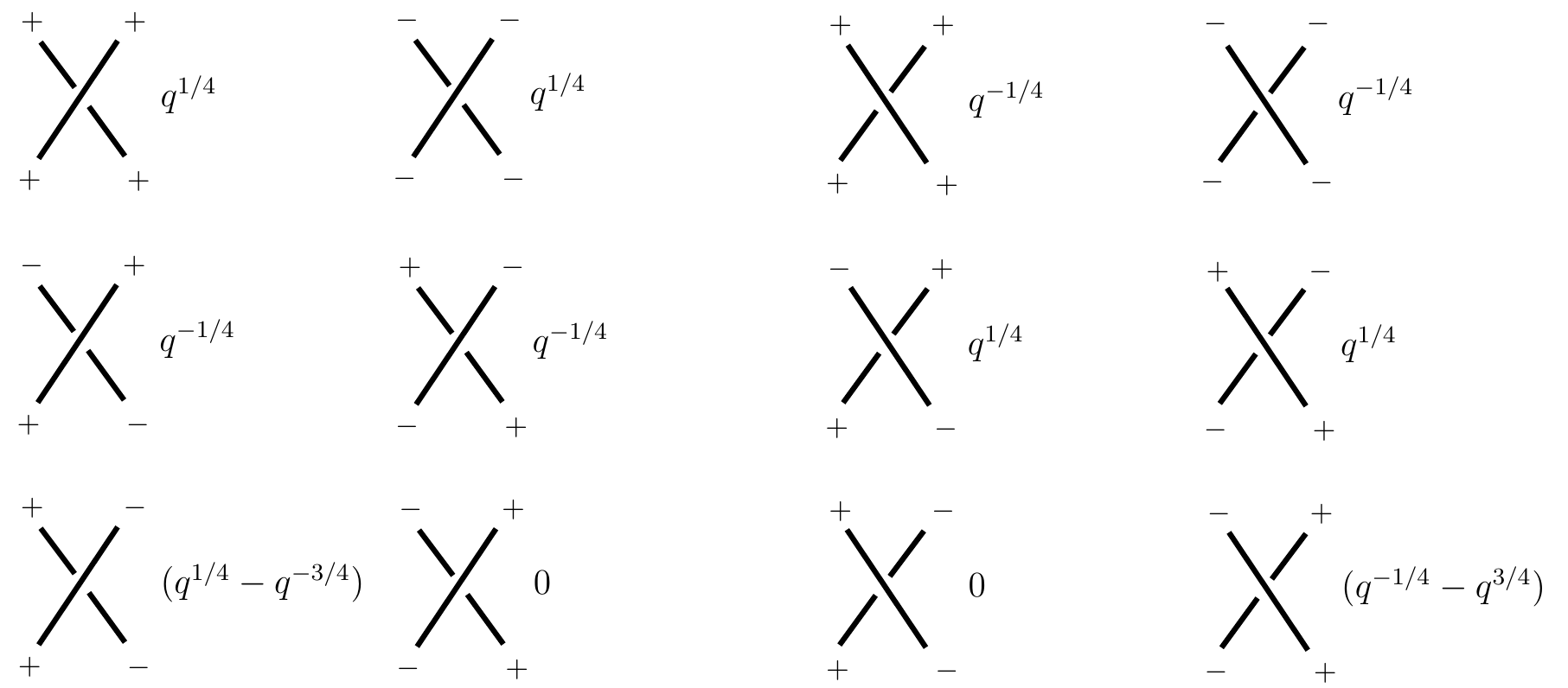}
 \end{center}
\caption{\small  The weights of the vertex model for a simple crossing of two strands.    (The weights for configurations not shown are 0.)  \label{weightone}}
\end{figure}

\begin{figure}
 \begin{center}
   \includegraphics[width=3in]{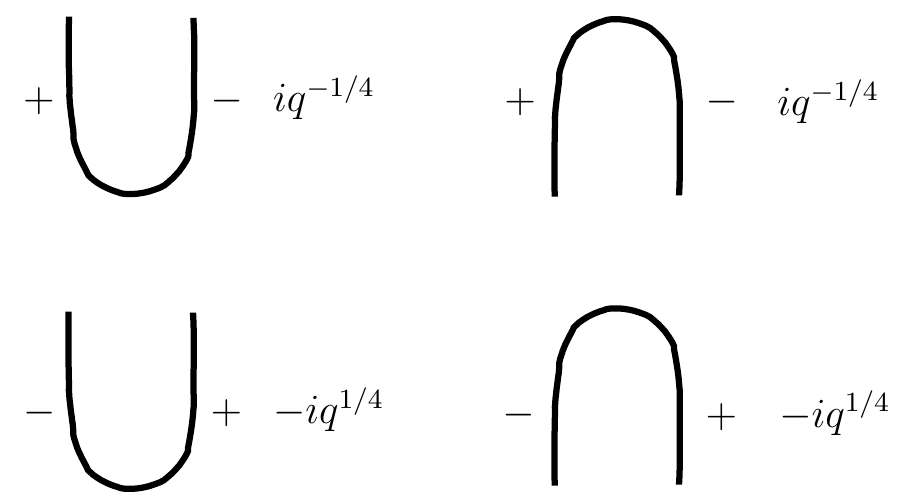}
 \end{center}
\caption{\small  The weights of the vertex model for a local maximum or minimum of the height function. (Weights not shown are 0.) \label{weighttwo}}
\end{figure}

For brevity, I will describe the ``vertex model'' (see \cite{Jonestwo} and also \cite{Kauffman}, p. 125).  One projects a knot to $\R^2$ in such
a way that the only singularities are  simple crossings
and so that the height function has only simple local maxima and minima  (fig. \ref{projection}).  One labels the intervals between crossings, maxima, and minima by a symbol $+$ or $-$.  One sums over all possible
labelings of the knot projection with simple weight functions given in figs. \ref{weightone} and \ref{weighttwo}.  The weights are functions of a variable $q$.  After summing
over all possible labelings  and weighting each labeling by the product of the weights attached to its crossings, maxima, and minima,
one arrives at a function of $q$.  The sum turns out to be an invariant of a framed knot.\footnote{A framing of a knot in $\R^3$  is
a trivialization of the normal bundle to the knot.  If a knot is given with a projection to $\R^2$, then the normal direction
to $\R^2$ gives a framing.   A change of framing multiplies the sum that comes from the vertex
model by an integer power of $q^{3/4}$.  A knot in $\R^3$ can be given
a canonical framing, and therefore the Jones polynomial can be expressed as an invariant of a knot without
a choice of framing.  But this is not always convenient.}       This invariant is a Laurent polynomial in $q$ (times a fixed fractional power of $q$
that depends on the framing).   It   is known as the Jones polynomial.

Clearly, given the rules stated in the figures, the Jones polynomial for a given knot is completely computable by a finite (but exponentially
long) algorithm.   The rules, however, seem to have come out of thin air.  Topological invariance is not obvious
and is proved by checking Reidemeister moves. 

\def\CS{{\mathrm{CS}}}
\def\Tr{{\mathrm{Tr}}}
\def\d{{\mathrm d}}

Other descriptions of the Jones polynomial were found during the same period, often involving mathematical physics.  The methods 
involved statistical mechanics,
braid group representations, quantum groups, two-dimensional conformal field theory, and more.  One notable fact was that conformal
field theory can be used \cite{TK} to generalize the constructions of Jones to the choice of an arbitrary simple 
Lie group\footnote{$G^\vee$ is a common notation for the Langlands or GNO dual of a simple Lie group $G$.  
Duality will later enter our story,
and we will have two descriptions involving a dual pair of groups $G$ and $G^\vee$.  We write $G^\vee$ for the group
that is used in the conformal field theory and Chern-Simons descriptions, because this will agree better with 
the usual terminology concerning the Langlands correspondence.  Similarly, we write $R^\vee$ for a representation of $G^\vee$ and $E^\vee$ for
a $G^\vee$ bundle.}  $G^\vee$ with a labeling of a knot
(or of each component of a link) by an irreducible representation $R^\vee$ of $G^\vee$.  The original 
Jones polynomial is the case that $G^\vee=SU(2)$ and $R^\vee$ is
the two-dimensional representation.  

With these and other clues, it turned out \cite{WittenJones} 
that the Jones polynomial can be described in
three-dimensional quantum gauge theory.   Here we start with a compact simple gauge group  $G^\vee$ 
(to avoid minor details, take $G^\vee$ to be connected
and simply-connected) and a trivial\footnote{If $G^\vee$ is connected and simply-connected, then inevitably any $G^\vee$-bundle over a three-manifold
is trivial.}  $G^\vee$-bundle $E^\vee\to W$, where $W$ is an oriented three-manifold.  Let $A$ be a connection on $E^\vee$.  The only gauge-invariant function of $A$
that we can write by integration over $W$ of some local expresssion, assuming no structure on $W$ except an orientation, is the Chern-Simons function
\begin{equation}\label{cs}\CS(A)=\frac{1}{4\pi}\int_W\Tr\left(A\wedge \d A+\frac{2}{3}A\wedge A\wedge A\right). \end{equation}
Even this function is only gauge-invariant modulo a certain fundamental period. 
In (\ref{cs}), $\Tr$ is an invariant and nondegenerate quadratic form on the Lie algebra of $G^\vee$, normalized so that $\CS(A)$ 
 is gauge-invariant mod $2\pi\Z$.   For $G^\vee=\SU(n)$ (for some $n\geq 2$), we can take $\Tr$ to be the trace in the $n$-dimensional representation.  

The Feynman path integral is now formally an integral over the infinite-dimensional space $U$ of connections:
\begin{equation}\label{delf} Z_k(W)=\frac{1}{\mathrm{vol}}\int_U DA\,\exp(ik \CS(A)). \end{equation}
This is a basic construction in quantum field theory, though unfortunately unfamiliar from a mathematical point of view.
Here
$k$ has to be an integer since $\CS(A)$ is only gauge-invariant modulo $2\pi\Z$.   $Z_k(W)$ is defined with no structure on $W$ except an orientation,
so it is an invariant of the oriented three-manifold $W$.  (Here and later, I ignore some details.  $W$ actually has to be ``framed,'' as one learns
if one follows the logic of ``renormalization theory.''  Also, formally $\mathrm{vol}$ is the volume of the infinite-dimensional group of gauge transformations.)

To include a knot -- that is an embedded oriented circle $K\subset W$ -- we make use of the {\it holonomy}
of the connection $A$ around $W$, which we denote  $\mathrm{Hol}(A,K)$. We pick
an irreducible  representation $R^\vee$ of $G^\vee$ and define
\begin{equation}\label{wilson}\W_{R^\vee}(K)=\Tr_{R^\vee}\,\mathrm{Hol}_K(A)=\Tr_{R^\vee} P\exp\left(-\oint_KA\right), \end{equation}
where the last formula is the way that physicists often denote the  trace of the holonomy.   In the context of quantum field theory, the trace of the holonomy
is usually called the Wilson loop operator.  Then we define a natural invariant of the pair $W,K$:
\begin{equation}\label{ilson}Z_k(W;K,R^\vee)=\frac{1}{\mathrm{vol}}\int_U DA \exp(ik\CS(A))\, \W_{R^\vee}(K). \end{equation}
(Again, framings are needed.)

If we take $G^\vee$ to be $\SU(2)$ and $R^\vee$ to be the two-dimensional representation, then $Z_k(W;K,R^\vee)$ turns out to be the Jones
polynomial, evaluated at\footnote{The analog of this for any simple  $G^\vee$  is $q=\exp(2\pi i/n_{\frak g}(k+h^\vee))$, where $h^\vee$ is the dual Coxeter number of $G^\vee$ and $n_{\frak g}$ is
the ratio of length squared of long and short roots of $G^\vee$.  In the dual description that we come to later, $q$ is always
the instanton-counting parameter.}  
\begin{equation}\label{fils}q=\exp\left(\frac{2\pi i}{k+2}\right). \end{equation}
  This statement is justified by making contact with two-dimensional
conformal field theory, via the results of \cite{TK}.  For a particularly direct way to establish the relation to the Knizhnik-Zamolodchikov equations of conformal field theory, see \cite{FK}.
 This relationship between three-dimensional gauge theory and two-dimensional
conformal field theory has also been important in condensed matter physics, in studies of the quantum Hall effect and related phenomena.

This approach has more or less the opposite virtues and drawbacks of standard approaches to the Jones polynomial.  No projection
to a plane is chosen, so topological invariance is obvious (modulo standard quantum field theory machinery), but it is not clear how
much one will be able to compute.  In other approaches, like the vertex model, there is an explicit finite algorithm for computation, but
topological invariance is obscure.

Despite the manifest topological invariance of this approach to the Jones polynomial, there were at least two things that many knot theorists
did not like about it.  One was simply that the framework of integration over function spaces -- though quite familiar to physicists -- is unfamiliar
mathematically.  (A version of this problem is one of the Clay Millennium Problems.)  The second is that this method does not give a clear
approach to understanding why the usual quantum knot invariants are Laurent polynomials in $q$ (and not just functions of an integer $k$).  
From some points of view, this is considered sufficiently important that it is part of the name ``Jones polynomial.''   Other approaches to the Jones polynomial -- such as the vertex model
that we started with -- do not obviously give a topological invariant but do obviously give a Laurent polynomial.

Actually, for most three-manifolds, the answer that comes from the gauge theory is the right one. It is special to knots in $\R^3$ that the natural
variable is $q=\exp(2\pi i/(k+2))$ rather than $k$.  The quantum knot invariants on a general three-manifold $W$ are naturally defined only for
an integer $k$ and do not  have natural analytic continuations to functions of\footnote{An analytic continuation can be made away from integer $k$, using ideas we explain
later.  On a generic three-manifold, the continued function has an essential singularity at $k=\infty$ with Stokes phenomena.  It is not a function of $q$.}  $q$.    This has been the traditional understanding: the gauge theory gives
directly a good understanding on a general three-manifold $W$, but if one wants to understand from three-dimensional gauge theory some of the special
things that happen for knots in $\R^3$, one has to begin by relating the gauge theory to one of the other approaches, for instance via conformal field theory.

However, a little over a decade ago, two developments gave clues that there should be another explanation.  
One of these developments was Khovanov homology, which will be the topic of the second lecture.  The other development, which started at roughly the
same time, was the ``volume conjecture'' \cite{Kash, MM,MMOTY,KT,Gu,Mm}.  What I will explain in this lecture started by trying to
understand the volume conjecture.  I should stress that I have not succeeded in finding a quantum field theory explanation for the volume conjecture.\footnote{I am not even entirely convinced that it is true.  What was found in \cite{Analytic} is that the volume conjecture for a certain knot is valid
if and only if a certain invariant of that knot is non-zero. (This invariant is the coefficient of a thimble associated to a flat $G^\vee_\C=SL(2,\C)$ connection
of maximal volume when the standard real integration cycle $U$ is expressed in terms of Lefschetz thimbles, along the lines of eqn. (\ref{nexth}).) It is not clear why this invariant is nonzero
for all knots.}
However, just understanding a few preliminaries concerning the volume conjecture led to a new point of view on the Jones polynomial.
This is what I aim to explain.   Since this is the case, I will actually not give a precise statement of the volume conjecture.

To orient ourselves, let us just ask how the basic integral
\begin{equation}\label{delff} Z_k(W)=\frac{1}{\mathrm{vol}}\int_U DA\,\exp(ik \CS(A))\end{equation}
behaves for large $k$.  It is an infinite-dimensional analog of a finite-dimensional oscillatory integral such as the one that defines the Airy function
\begin{equation}\label{elf}F(k;t)=\int_{-\infty}^\infty \d x\,\exp(ik(x^3+tx)), \end{equation}
where we assume that $k$ and $t$ are real.  Taking $k\to\infty$ with fixed $t$, the integral vanishes exponentially fast if there are no real
critical points ($t>0$) and is a sum of oscillatory contributions of real critical points if there are any ($t<0$).  The same logic applies to the infinite-dimensional
integral for $Z_k(W)$.   The critical points of $\CS(A)$ are flat connections, corresponding to homomorphisms $\rho:\pi_1(W)\to G$, so the
asymptotic behavior of $Z_k(W)$ for large $k$ is given by a sum of oscillatory contributions associated to such homomorphisms. 
(This has been shown explicitly in examples \cite{FG,LJe}.)  

The volume conjecture arises if we specialize to knots in $\R^3$, so that $k$ does not have to be an integer.  Usually the case $G^\vee=\SU(2)$ is
assumed and we let $R^\vee$ be the $n$-dimensional representation of $\SU(2)$.  (The corresponding knot invariant is called the colored Jones polynomial.)
Then we take $k\to\infty$ through noninteger values, with fixed $k/n$.   A choice that is sufficient for our purposes is to set $k=k_0+n$, where $k_0$ is 
a fixed complex number and we take $n\to\infty$ (through integer values).  The large $n$ behavior is now a sum of contributions from {\it complex}
critical points.  By a complex critical point, I mean simply a critical point of the analytic continuation of the function $\CS(A)$.

We make this analytic continuation simply by replacing the Lie group $G^\vee$ with its complexification $G^\vee_\C$, replacing the $G^\vee$-bundle $E^\vee\to W$
with its complexification, which is a $G^\vee_\C$ bundle $E^\vee_\C\to W$, and replacing the connection $A$ on $E^\vee$ by a connection $\A$ on $E^\vee_\C$,
which we can think of as a complex-valued connection.  Once we do this, the function $\CS(A)$ on the space $U$ of connnections on $E^\vee$ can
be analytically continued to a holomorphic function $\CS(\A)$ on $\U$, the space of connections on $E^\vee_\C$.  This function is defined by the
``same formula'' with $A$ replaced by $\A$:
\begin{equation}\label{csc}\CS(\A)=\frac{1}{4\pi}\int_W\Tr\left(\A\wedge \d \A+\frac{2}{3}\A\wedge \A\wedge \A\right). \end{equation}
On a general three-manifold $W$, a critical point of $\CS(\A)$ is simply a complex-valued flat connection, corresponding to a homomorphism
$\rho:\pi_1(W)\to G^\vee_\C$.

In the case of the volume conjecture with $W=\R^3$, the fundamental group is trivial, but we are supposed to also include a holonomy or Wilson
loop operator $\W_{R^\vee}(K)=\Tr_{R^\vee}\,\mathrm{Hol}_K(A)$, where $R^\vee$ is the $n$-dimensional representation of $SU(2)$. When we take $k\to\infty$
with fixed $k/n$, this holonomy factor affects what we should mean by a critical point.\footnote{\label{ifinstead} If instead we take $k\to\infty$
with fixed $n$, we do not include $W_{R^\vee}(K)$ in the definition of a critical point; we simply view it as a function that can be evaluated at a critical
point of $\CS(\A)$.  We will follow this second procedure later.}   A full explanation would take us too far afield, and instead I will just say the answer: the right notion of a complex critical
point for the colored Jones polynomial is a homomorphism $\rho:\pi_1(W\backslash K)\to G^\vee_\C$, with a monodromy around $K$ whose conjugacy
class is determined by the ratio $n/k$.  What is found in work on the ``volume conjecture'' is that typically the colored Jones polynomial for $k\to\infty$
with fixed $n/k$ is determined by such a complex critical point.  

Physicists know about various situations (involving ``tunneling'' problems) in which a path
integral is dominated by a complex critical point, but usually this is a complex critical point that
makes an exponentially small contribution.  There is a simple reason for this.  Usually in
quantum mechanics one is computing a probability amplitude.  Since probabilities cannot be bigger than
1, the contribution of a complex critical point to a probability amplitude can be exponentially small
but it cannot be exponentially large.  What really surprised me about the volume conjecture
is that, for many knots (knots with hyperbolic complement, in particular), the dominant critical
point makes an exponentially {\it large} contribution.  In other words, the colored Jones polynomial
is a sum of oscillatory terms  for $n\to\infty$, $k=k_0+n$ if $k_0$ is an integer, but it grows
exponentially in this limit as soon as $k_0$ is not an integer.  (Concretely, that is because 
$k\CS(\A)$ evaluated at the appropriate critical point has a negative imaginary part, so
$\exp(ik\CS(\A))$ grows exponentially for large $k$.)

There is no contradiction with the statement that quantum mechanical probability amplitudes cannot
be exponentially large, because as soon as $k_0$ is not an integer, we are no longer studying a physically
sensible quantum mechanical system. But it seemed puzzling that making $k_0$ non-integral,
even if  still real, can change the large $n$ behavior so markedly.  However, it turns out that a
simple one-dimensional integral can do the same thing:
\begin{equation}\label{ix} I(k,n)=\int_0^{2\pi} \frac{\d\theta}{2\pi} e^{ik\theta}e^{2in\sin\theta}. \end{equation}
We want to think of $k$ and $n$ as analogs of the integer-valued parameters in Chern-Simons
gauge theory that we call by the same names.  (In our model problem, $k$ is naturally an integer,
but there is no good reason for $n$ to be an integer.  So the analogy is not perfect.)
If one takes $k,n$ to infinity with a fixed (real) ratio and maintaining integrality of $k$, the integral $I(k,n)$ has an oscillatory behavior,
dominated by the critical points of the exponent $f=k\theta+2n\sin\theta$, if $k/n$ is such that
there are critical points for real $\theta$.  Otherwise, the integral vanishes exponentially fast for large $k$.

Now to imitate the situation considered in the volume conjecture, we want to analytically continue
away from integer values of $k$.
The  integral $I(k,n)$ obeys Bessel's equation (as a function of $n$) for any integer $k$.
We want to think of Bessel's equation as the analog of the ``Ward identities'' of quantum field
theory, so in the analytic continuation of $I(k,n)$ away from integer $k$, 
we want to preserve Bessel's equation.
The proof of Bessel's equation involves integration by parts, so it is important that we are integrating
all the way around the circle and that the integrand is continuous and single-valued on the circle.   That is
why $k$ has to be an integer.

\begin{figure}
 \begin{center}
   \includegraphics[width=3in]{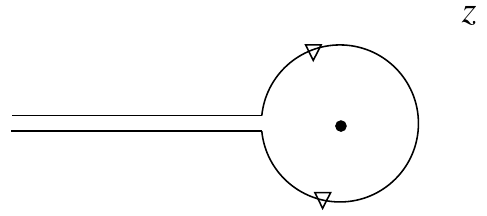}
 \end{center}
\caption{\small   The contour used in analytic continuation of the solution of Bessel's equation. \label{newcontour}}
\end{figure}

The analytic continuation of $I(k,n)$, preserving Bessel's equation, was known in the 19th century.
We first set $z=e^{i\theta}$, so our integral becomes
\begin{equation}\label{morx}I(k,n)=\oint\frac{\d z}{2\pi i}z^{k-1}\exp\left(n(z-z^{-1})\right). \end{equation}
Here the integral is over the unit circle in the $z$-plane.  At this point, $k$ is still an integer.  We want
to get away from integer values while still satisfying Bessel's equation.  If $\mathrm{Re}\,n>0$, this can be done
by switching to the integration cycle shown in fig. \ref{newcontour}.
 
The integral on the new cycle converges (if $\mathrm{Re}\,n>0$), and it agrees with the original integral on the circle if $k$ is an integer, since
the extra parts of the cycle cancel.  But the new cycle gives a continuation away from integer $k$, still obeying Bessel's equation.  
There is no difficulty in the integration by parts used to prove Bessel's equation, since the integral on the chosen cycle is rapidly convergent at infinity.  

How does the integral on the new cycle behave in the limit $k,n\to\infty$ with fixed $k/n$?  If $k$ is an integer and $n$ is real, the integral is
oscillatory or exponentially damped, as I have stated before, depending on the  ratio $k/n$.  But as soon as $k$ is not an integer (even if $k$ and $n$
remain real), the large $k$ behavior with fixed $k/n$ can grow exponentially, for a certain range of $k/n$, rather as one finds for the colored Jones
polynomial. Unfortunately, even though it is elementary, to fully explain this statement would involve a bit of a digression.  (Details can be found,
for example, in \cite{Analytic}, section 3.5.)   Here I will just explain the technique that one can use to make this analysis, since this will show
the technique that we will follow in taking a new look at the Jones polynomial.  

We are trying to do an integral of the generic form
\begin{equation}\label{omely}\int_\Gamma \frac{\d z}{2\pi i z}\exp(kF(z)),\end{equation}
where $F(z)$ is a holomorphic function and $\Gamma$ is a cycle, possibly not compact, on which the integral
converges.   In our case, 
\begin{equation}\label{colum}F(z)=\log z +\lambda(z-z^{-1}),~~\lambda=n/k. \end{equation}
We note that because of the logarithm, $F(z)$ is multivalued.  To make the analysis properly, we should work on a cover of the punctured
$z$-plane parametrized by $w=\log z$ on which $F$ is single-valued:
\begin{equation}\label{melx} F(w)=w+\lambda(e^w-e^{-w}). \end{equation}
The next step is to find  a useful description of all possible cycles on which the desired integral, which now is
\begin{equation}\label{elx}\int_\Gamma\frac{\d w}{2\pi i}\exp(kF(w)), \end{equation}
converges.  

Morse theory gives an answer to this question.  We consider the function $h(w,\bar w)=\mathrm{Re}(k F(w))$ as a Morse function.  Its critical points
are simply the critical points of the holomorphic function $F$ and so in our example they obey
\begin{equation}\label{mixo}1+\lambda(e^w+e^{-w})=0. \end{equation}
The key step is now the following.  To every critical point $p$ of $F$, we can define an integration cycle $\Gamma_p$, called a Lefschetz
thimble, on which the integral we are trying to do converges.   Moreover, the $\Gamma_p$ give a basis of integration cycles on
which this integral converges, since  they give a basis of the homology of the $w$-plane relative
to the region with $h\to -\infty$.  
(We assume that the critical points of $F$ are all non-degenerate, as is the case in our
example.  Also, we  assume that $F$ is  sufficiently generic so that the equation (\ref{zinc}) introduced momentarily has no solutions
interpolating from one critical point at $t=-\infty$ to another at $t=+\infty$.  If $F$ varies as a function of some parameters,
then in real codimension 1, such interpolating solutions do appear; there is then a Stokes phenomenon -- a jumping in the basis of the relative homology given by the $\Gamma_p$.)

 \begin{figure}
 \begin{center}
   \includegraphics[width=1.7in]{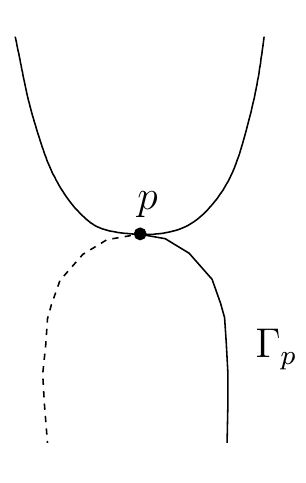}
 \end{center}
\caption{\small   The Lefschetz thimble associated to a critical point.  The critical point is a saddle point and the thimble is the union of downward flows
that start at this saddle.\label{thimble}}
\end{figure}

In fact, since $h$ is the real part of a holomorphic function, its critical points are all saddle points, not local maxima or minima.  The Lefschetz
thimble associated to a given critical point $p$ is defined by ``flowing down'' from $p$ (fig. \ref{thimble}), via the gradient flow equation of Morse
theory.  We could use any complete Kahler metric on the $w$-plane in defining this equation, but we may as well use the obvious flat metric
$\d s^2=|\d w|^2$.  The gradient flow equation is then 
\begin{equation}\label{zinc}\frac{\d w}{\d t}=-\frac{\partial h}{\partial\bar w}, \end{equation}
where $t$ is a new ``time'' coordinate.  The Lefschetz thimble $\Gamma_p$ associated to a critical point $p$ is defined as the space of all
values at $t=0$ of solutions of the flow equation on the semi-infinite interval $(-\infty,0]$ that start at $p$ at $t=-\infty$.  For example, $p$ itself
is contained in $\Gamma_p$, because it is the value at $t=0$ of the trivial solution of the flow equation that is equal to $p$ for all times.  
A non-constant solution that approaches $p$ for $t\to-\infty$ is exponentially close to $p$ for large negative $t$.  The coefficient of the exponentially
small term  in a particular solution determines how far the flow reaches by time $t=0$ and therefore what point on $\Gamma_p$ is represented
by this particular flow.  

$\Gamma_p$ is not compact, but the integral 
\begin{equation}\label{mert}I_p=\int_{\Gamma_p}\frac{\d w}{2\pi i}\exp(k F(w)) \end{equation}
converges, since $h=\mathrm{Re}(kF(w))$ goes to $-\infty$ at infinity along $\Gamma_p$.  Moreover,  when restricted to $\Gamma_p$, $h$ has a unique maximum, which is at the point $p$.   This statement leads to a straightforward answer for the large $k$ behavior of the integral $I_p$:
\begin{equation}\label{sufuf}I_p\sim \exp(k F(p))\left(c_0k^{-1/2}+c_1k^{-3/2}+\dots\right), \end{equation}
where the coefficients $c_0,c_1,\dots$ in the asymptotic expansion can be computed by classical methods.

Any other cycle $\Gamma$ on which the integral converges can be expanded as a linear combination of the Lefschetz thimbles:
\begin{equation}\label{expth}\Gamma=\sum_p a_p \Gamma_p, ~~a_p\in \Z. \end{equation}
After computing the integers $a_p$, it is straightforward to determine the large $k$ asymptotics of an integral
\begin{equation} \label{nexth}\int_\Gamma\frac{\d w}{2\pi i}\exp(kF(w)). \end{equation}
It is simply given by the contributions of those critical points $p$ for which $h(p)$ is maximal under the condition that $a_p\not=0$.
Applying this procedure to our example related to the Bessel function, we get the answer that I claimed before: this integral has an asymptotic behavior similar to that
of the colored Jones polynomial.  The limit $n\to\infty$, $k=k_0+n$ is quite different depending on whether $k_0$ is an integer.  (Concretely, if $k_0$ is not an integer,
the large $n$ behavior is dominated by two Lefschetz thimbles whose contributions cancel if $k_0$ is an integer.) 

At this stage, I hope it is fairly clear what we should do to understand the analytic continuation to non-integer $k$ of the quantum invariants of knots in $\R^3$,
and also to understand the asymptotic behavior of the colored Jones polynomial that is related to the volume conjecture.  We should define Lefschetz thimbles in the space
$\U$ of complex-valued connections, or more precisely in a cover of this space on which $\CS(\A)$ is single-valued,
and in the gauge theory definition of the Jones polynomial, we should replace the integral over the space $U$ of real connections
with a sum of integrals over Lefschetz thimbles.  

\def\Re{{\mathrm{Re}}}
\def\Im{{\mathrm{Im}}}
However, it probably is not clear that this will actually lead to a useful new viewpoint on the Jones polynomial.  This depends on a few additional facts. 
To define the Lefschetz thimbles that we want, we need to consider a gradient flow equation on the infinite-dimensional space $\U$ of complex-valued connections,
with $\Re(ik\CS(\A))$ as a Morse function.\footnote{For a more detailed explanation of the following, see \cite{Wittengauge}.}  
Actually, I want to first practice with the case of gradient flow on the infnite-dimensional space $U$ of real connections (on a $G^\vee$-bundle
$E^\vee\to W$, $W$ being a three-manifold)
with the Morse function being the real Chern-Simons function $\CS(A)$.  This case is important in Floer theory of three-manifolds and in Donaldson theory of smooth
four-manifolds, so it is relatively familiar.  A Riemannian metric on $W$ induces a Riemannian metric on $U$ by 
\begin{equation}\label{thorg} |\delta A|^2=-\int_W\Tr\,\delta A\wedge \star \delta A, \end{equation}
where $\star=\star_3$ is the Hodge star operator acting on differential forms on $W$.    
We will use this metric in defining a gradient flow equation on $U$, with Morse function $\CS(A)$.

The flow equation will be a differential equation on a four-manifold $M=W\times \R$, where $\R$ is parameterized by the ``time''; one can think of the flow as evolving a three-dimensional
connection in ``time.''  Concretely, the flow equation is 
\begin{equation}\label{tholf}\frac{\partial A}{\partial t}=-\frac{\delta \CS(A)}{\delta A}=-\star_3 F, \end{equation}
where $F=\d A+A\wedge A$ is the curvature.  Now a couple of miracles happen.  This equation has {\it a priori} no reason to be elliptic or to have four-dimensional
symmetry. But it turns out that the equation is actually a gauge-fixed version of the instanton equation $F^+=0$, which is elliptic modulo the gauge group and has the full
four-dimensional symmetry (that is, it is naturally-defined on any oriented Riemannian  four-manifold $M$, not necessarily of the form $W\times \R$ for some $W$).  
These miracles are well-known to researchers on Donaldson and Floer theory, where they play an important role.

\def\ad{{\mathrm{ad}}}
It turns out that similar miracles happen in gradient flow on the space $\U$ of complex-valued connections, endowed with the obvious flat Kahler metric
\begin{equation}\label{worf}|\delta \A|^2=-\int_W\Tr\,\delta \A\wedge\star\delta \bar\A. \end{equation}
This equation is a gauge-fixed version (with also the moment map set to 0, in a sense explained in \cite{Wittengauge}) of an elliptic differential equation
that has full four-dimensional symmetry.  This equation can be seen as a four-dimensional cousin of Hitchin's celebrated equation in two dimensions.
It is an equation for a pair $A,\phi$, where $A$ is a real connection on a $G^\vee$-bundle $E^\vee\to M$, $M$ being an oriented four-manifold, and $\phi$ is a one-form on $M$
with values in $\ad(E^\vee)$.  The equations (for simplicity I take $k$ real) are 
\begin{equation}\label{delft}F-\phi\wedge\phi=\star\d_A\phi, ~~~\d_A\star\phi = 0. \end{equation}
They can be viewed as flow equations for the complex-valued connection $\A=A+i\phi$ on the three-manifold $W$. 

There is a happy coincidence; these equations, which sometimes have been called the KW equations,
 arise in a certain twisted version of maximally supersymmetric Yang-Mills theory ($\N=4$ super Yang-Mills theory) in four dimensions \cite{KW}.
We will see shortly why this relationship is relevant.  For recent mathematical work on these equations, see \cite{Taubes,Taubestwo} and also \cite{GU}.

Now we can define a Lefschetz thimble for any choice
 of a complex flat connection $\A_\rho$ on $M$, associated to a homomorphism
 $\rho:\pi_1(M)\to G^\vee_\C$.  We work on the four-manifold $M=W\times \R_+$,
 where $\R_+$ is the half-line $t\geq 0$, and define the thimble $\Gamma_\rho$ to consist of all complex  connections $\A=A+i\phi$ that are boundary
 values (at the finite boundary of $M$ at $W\times\{t=0\}$) of solutions on the KW equations on $M$ that approach $\A_\rho$ at infinity.
 
\begin{figure}
 \begin{center}
   \includegraphics[width=3in]{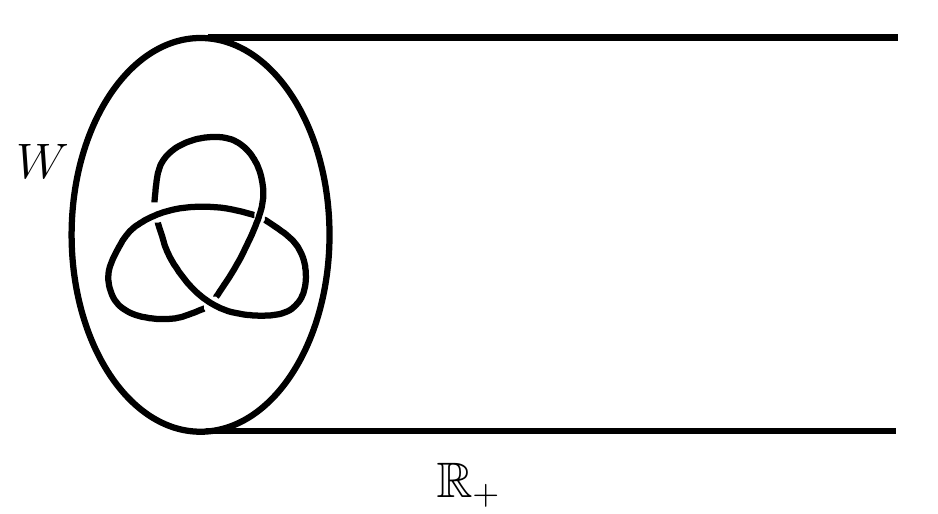}
 \end{center}
\caption{\small   A four-manifold $M=W\times \R_+$, with a knot $K$ embedded in its boundary $W\times\{0\}$. \label{longone}}
\end{figure}

 For a general $M$, there are various choices of $\rho$ and some rather interesting issues that have not yet been unraveled.  But now we 
 can see what is special about knots in $\R^3$.  Since the fundamental group of $\R^3$ is trivial,\footnote{\label{loco} As in footnote \ref{ifinstead}, there are two fruitful approaches to the present subject, which differ by whether in getting a semiclassical limit,
 we {\it (i)} take $k\to \infty$ with fixed $n/k$, or {\it (ii)} take $k\to\infty$ with fixed $n$.  In approach {\it (i)}, a complex critical point is
 a complex flat connection on the knot complement $W\backslash K$.  We followed this approach in our discussion of the volume conjecture.
 In approach {\it (ii)}, which is more convenient here and in Lecture Two below, a complex critical point is simply a complex flat connection
 on $W$.  In approach {\it (ii)}, the holonomy function $\W_{R^\vee}(K)$ does not enter the definition of a Lefschetz thimble, but is just one factor
 in the function that we want to integrate over the thimble.  
   The two
viewpoints are described more thoroughly in \cite{WittenK}.}  any complex flat connection on $\R^3$ is
 equivalent to the trivial one, $\A=0$.  Hence there is only one Lefschetz thimble $\Gamma_0$, and any integration cycle is a 
multiple of this one.  So instead of integration over 
 $U$ to define the Jones polynomial, we can define the quantum knot invariants  by integration over $\Gamma_0$:
 \begin{equation}\label{turkey} Z_k(\R^3;K,R^\vee)=\frac{1}{\mathrm{vol}}\int_{\Gamma_0}D\A \,\exp\left(ik\CS(
\A)\right)\cdot \W_{R^\vee}(K). 
 \end{equation} Here the holonomy function $\W_{R^\vee}(K)$ is viewed as a function on the Lefschetz thimble; in other words,
 it is evaluated for the connection $\A=A+i\phi$ restricted to $W\times \{0\}$ (so in fig. \ref{longone}, the knot $K$ is placed
 on the boundary of $W\times \R_+$; as in footnote \ref{loco}, $K$ does not enter the definition of the Lefschetz thimble).  
 This definition explains why the quantum invariants of knots in $\R^3$ can be analytically continued away from roots of unity.  Indeed,
 the function $\CS(\A)$ is well-defined and single-valued on the Lefschetz thimble $\Gamma_0$, so there is no 
 reason to restrict to the case that $k$
 is an integer.

 The formula (\ref{turkey}) for the Jones polynomial and its cousins may appear to be purely formal, but there is a reason that we can say something
 about it.  As I have already observed, the KW equations arise in  $\N=4$ super Yang-Mills theory in four dimensions.  This
 theory has a ``twisted'' version that localizes on the space of solutions of the KW equations. The space of all such solutions on $M=\R^3\times \R_+$,
 with the requirement that $\A\to 0$ at $\infty$, is simply our Lefschetz thimble $\Gamma_0$.  The upshot is that the quantum  invariants of
 a knot in $\R^3$ can be computed from a path integral of $\N=4$ super Yang-Mills theory in four dimensions, with a slightly subtle boundary
 condition \cite{GW} along the boundary at $\R^3\times \{0\}$.
 
 This is not yet obviously useful, but one more step brings us into a more accessible world, and also gives a new explanation
 of why the quantum knot invariants are Laurent polynomials in the variable $q$.  The step in question was 
 also a key step in \cite{KW} and more
 generally in most of the work of physicists on the supersymmetric gauge theory in question.  This is electric-magnetic 
 duality, the four-dimensional
 analog of mirror symmetry in two dimensions.  $\N=4$ supersymmetric Yang-Mills theory with gauge group $G^\vee$ and 
 ``coupling parameter''
 $\tau$ is equivalent to the same theory with $G^\vee$ replaced by its Langlands or Goddard-Nuyts-Olive dual group, which we simply call  $G$, and coupling parameter $\tau^\vee$ replaced by $\tau
 =-1/n_{\frak g}\tau^\vee$ (here $n_{\frak g}$
 is the ratio of length squared of long and short roots of $G$ or equivalently of $G^\vee$).  
 
 To find a dual description in our problem, we need to ask what happens  under
 the duality to the boundary condition at $\R^3\times \{0\}$.  (The analog of this question in mirror symmetry may be more familiar: what Lagrangian submanifold is mirror to a given
 coherent sheaf?)  For the boundary condition that is related to the Lefschetz thimble, the dual boundary condition was described
 some years ago
 in \cite{GW}.  It is somewhat unusual and will be described in the next lecture.    For now, I will just
 say that this boundary condition has the formal properties of a standard local elliptic boundary condition and has the effect of
 reducing to finite-dimensional spaces of solutions of the KW equations. 
 
 In the situation of fig. \ref{longone}, after making the duality transformation, the moduli space of solutions has expected dimension
 0 and to evaluate  $Z_k(\R^3;K,R^\vee)$, we just have to ``count''  (with signs, as in Donaldson theory) the number $b_n$
 of solutions for a given value $n$ of the instanton number (for $G=\SU(n)$, the instanton number is the second Chern class).
 The boundary conditions depend on the knot $K$ and on the representation $R^\vee$ by which it is labeled.  This is the only
 way that $K$ and $R^\vee$ enter in this dual description.  The path integral gives
 \begin{equation}\label{morfo} Z(q;K,R^\vee)=\sum_n b_n q^n, \end{equation}
 where $q$ was defined in eqn. (\ref{fils}).  This exhibits the Jones polynomial
and the related quantum invariants of knots in three dimensions as
``Laurent polynomials'' in $q$ with integer coefficients. I put
``Laurent polynomials'' in quotes because the powers of $q$ are
shifted from integers in a way that depends only on the
representations, so for instance the Jones polynomial of a knot
with this normalization is $q^{1/2}$ times a Laurent polynomial in $q$.

 I have changed
 notation slightly to write the knot invariant as $Z(q;K,R^\vee)$ rather than $Z_k(\R^3;K,R^\vee)$, since this formula only works
 in this simple way on the three-manifold $\R^3$, and also in this description the natural variable is $q$ rather than $k$.

 The formula (\ref{morfo}), in which $q$ appears as an instanton-counting parameter, can be viewed as a response to a
 challenge raised in \cite{AtiyahWeyl}, p. 299.  The challenge was to find a description of the Jones polynomial in which $q$ would
 be associated to instanton number (the integral of the second Chern class)  in four dimensions.  
 
 \section{Lecture Two}

 As we discussed in the last lecture, quantum knot invariants of a simple Lie group $G^\vee$ on a three-manifold $W$ can
 be computed by counting solutions of a certain system of elliptic partial differential equations, with gauge group the dual
 group $G$,  on the four-manifold $M=W\times \R_+$.  The equations are the KW equations
 \begin{equation}\label{melf} F-\phi\wedge\phi=\star\d_A\phi,~~~\d_A\star\phi=0 \end{equation}
 where $A$ is a connection on a $G$-bundle $E\to M$ and $\phi\in \Omega^1(M,\ad(E))$.  The boundary conditions
 at the finite end of $M=W\times \R_+$ depend on the knot, as indicated in fig. \ref{longone}.  The boundary conditions\
 at the infinite end of $M$ say that $\A=A+i\phi$ must approach a complex-valued flat connection.  Exactly what we have to do
 depends on what we want to get, but in one very important case there is a simple answer.  For $W=\R^3$, meaning
 that we are studying knots in $\R^3$, a flat connection is gauge-equivalent to zero and we require that $\A\to 0$ at infinity,
 in other words $A,\phi\to 0$.  In this lecture, we are only going to discuss the case that $W=\R^3$.
 
 For $W=\R^3$, the difference between $G^\vee$ and $G$ is going to be important primarily when they have different Lie
 algebras, since for instance there is no second Stieffel-Whitney class to distinguish $\SO(3)$ from $\SU(2)$.  So the
 difference will be most important if $G^\vee=\SO(2n+1)$ and $G=\Sp(2n+1)$, or vice-versa.   In fact, we will see later
 that something very interesting happens precisely for $G^\vee=\SO(2n+1)$ (or its double cover $\Spin(2n+1)$).
 
 To compute quantum knot invariants, we are supposed to ``count'' the solutions of the KW equations with fixed instanton
 number.  The instanton number is defined as
 \begin{equation}\label{murm} P=\frac{1}{8\pi^2}\int_M\Tr \, F\wedge F,\end{equation}
where the trace is an invariant quadratic form defined so that (for simply-connected $G$ on a compact
four-manifold $M$ without boundary), $P$ takes integer values.
For $G=\SU(n)$, we can take $\Tr$ to be the trace in the $n$-dimensional representation; then $P$ is the second Chern class. 
Just as in Donaldson theory, the ``count'' of solutions is made with signs.  The sign with which a given solution contributes is
the sign of the determinant of the linear elliptic operator that arises by linearizing the KW equations about a given solution. (For
physicists, this is the fermion determinant.)

Let $b_n$ be the ``number'' of solutions of instanton number $P=n$.  One forms the series
\begin{equation}\label{um}Z(q)=\sum_n b_nq^n. \end{equation}
One expects that $b_n$ vanishes for all but finitely many $n$.  Given this, $Z(q)$ (which depends on the knot $K$ and
a representation $R^\vee$, though we now omit these in the notation) is a Laurent polynomial in $q$ (times $q^c$
for some fixed  $c\in \Bbb Q$, as explained shortly).  For example, if $G^\vee=\SU(2)$ and the knot is labeled by the two-dimensional
representation, then $Z(q)$ is expected to be the Jones polynomial.

In all of this, the knot and representation are encoded entirely in the boundary condition at the finite end of $M$,
as sketched in fig. \ref{longone}.  The instanton number $P$ is an integer if $M$ is compact and without boundary,
but we are not in that situation.  To make $P$ into a topological invariant, we need a trivialization of the bundle $E$ at both the
finite and infinite ends of $M$.  The trivialization at the infinite end comes from the requirement that $A,\phi\to 0$ at infinity.
The trivialization at the finite end depends on the boundary condition, which I have not yet described.  With this boundary condition,
$P$ is offset from being integer-valued by a constant that only depends on the knots $K_i$ in $W$ and the representations 
$R_i^\vee$ labeling them.  This is why $Z(q)$ is not quite a Laurent polynomial in $q$, but is $q^c$ times such a Laurent
polynomial, where $c$ is completely determined by the representations $R_i^\vee$ and the framings of the $K_i$.

Given this description of the Jones polynomial and related knot invariants, I want to explain how to associate these knot invariants
with a homology theory (which is expected to coincide with Khovanov homology).
I should say that the original work by physicists associating vector spaces to knots was by Ooguri and Vafa \cite{OV}
(following earlier work associating vector spaces to homology cycles in a Calabi-Yau manifold
\cite{GV,GV2}).  After the invention of Khovanov homology of knots \cite{Khovanov}, a relation of the Ooguri-Vafa
construction to Khovanov homology was proposed \cite{GSV}.  What I will be summarizing here is
a parallel construction \cite{WittenK} in gauge theory language.
The arguments are probably more self-contained, though it is hard
to make this entirely clear in these lectures; the construction is
more uniform for all groups  and representations; and I believe
that the output is something that mathematicians will be able to
grapple with even without a full understanding of the underlying
quantum field theory. I should also say that my proposal for
Khovanov homology is qualititatively similar to ideas by Seidel and
Smith, Kronheimer and Mrowka, and probably others, and is expected to be mirror to  
a construction of Cautis and Kamnitzer.  See \cite{CK,SS,Kam,KM,KM2}
for references to this mathematical work.

\def\S{{\mathcal S}}
\def\V{{\mathcal V}}
Let $\S$ be the set of solutions of the KW equations.  (It is expected that for a generic embedding of a knot or link in $\R^3$,
the KW equations have only finitely many solutions and these are nondegenerate: the linearized operator has trivial kernel
and cokernel.)  We define a vector space $\V$ by declaring that for every $i\in\S$, there is a corresponding basis vector 
$|i\rangle$.  On $\V$, we will have two ``conserved quantum numbers,'' which will be the instanton number $P$ and a second
quantity that I will call the ``fermion number'' $F$.    I have already explained that $P$ takes values in $\Z+c$, where $c$ is
a fixed constant that depends only on the choices of representations and framings.  $F$ (which will be defined as a certain
Morse index) will be integer-valued.
The states $|i\rangle$ corresponding to solutions will be eigenstates of $F$.  We consider the state $|i\rangle$ to be
``bosonic'' or ``fermionic'' depending on whether it has an even or odd value of $F$.  So the operator distinguishing
bosons from fermions is $(-1)^F$.  $F$ will be defined so that if the solution $i$ contributed $+1$ to the counting of KW solutions,
then $|i\rangle$ has even $F$, and if it contributed $-1$, then $|i\rangle$ has odd $F$.  

Let us see how we would rewrite in this language the quantum knot invariant
\begin{equation}\label{urky}Z(q)=\sum_nb_nq^n. \end{equation}
Here a solution $i\in\S$ with instanton number $n_i$ and fermion number $f_i$ contributes $(-1)^{f_i}$ to $b_{n_i}$,
so it contributes $(-1)^{f_i}q^{n_i}$ to $Z(q)$.  So an equivalent formula is
\begin{equation}\label{rky} Z(q)=\sum_{i\in\S}(-1)^{f_i}q^{n_i}=\Tr_\V\,(-1)^Fq^P. \end{equation}

\def\H{{\mathcal H}}
So far, we have not really done anything except to shift things around.  However, on $\V$ we will also have a ``differential''
$Q$, which is an operator that commutes with $P$ but increases $F$ by 1, and obeys $Q^2=0$.   These statements mean
that we can define the {\it cohomology} of $Q$, which we denote as $\H$, and moreover that $\H$ is $\Z\times \Z$-graded, 
with the two gradings determined by  $P$ and $F$ (we simplify slightly, ignoring that the eigenvalues
of $P$ are really in  a coset $\Z+c$).

The importance of passing from $\V$ to $\H$ is that $\H$ is a topological invariant, while $\V$ is not.  If one deforms
a knot embedded in $\R^3$, solutions of the KW equations on $M=\R^3\times \R_+$ will appear and disappear, so $\mathcal V$
will change.  But $\H$ does not change.   This $\H$ is the candidate for the Khovanov homology.

Instead of defining $Z(q)$ as a trace in $\V$ via (\ref{rky}), we can define it as a trace in $\H$:
\begin{equation}\label{ky} Z(q)=\Tr_\H\,(-1)^Fq^P. \end{equation}
So here, $Z(q)$ is expressed as an ``Euler characteristic,'' i.e. as a trace in which bosonic and fermionic states cancel, in the
invariantly defined cohomology $\H$.  The reason that we can equally well write $Z(q)$ as a trace in $\V$ or in $\H$ is standard: the difference
 between $\V$ and $\H$ is that
in passing from $\V$ to $\H$, pairs of states disappear that make vanishing contributions to $Z(q)$.  (Such a pair consists of a bosonic state and a fermionic state with
the same value of $P$ and values of $F$ differing by 1.)

Defining the $\Z\times \Z$-graded vector space $\H$ and not just the trace $Z(q)$  adds information for two reasons.
One reason is simply that the fermion number $F$ is really $\Z$-valued, and this is part of the $\Z\times\Z$ grading of $\H$.
When we pass from $\H$ to the trace $Z(q)$, we only remember $F$ modulo 2, and here we lose some information.  The
second reason is that one can define natural operators acting on $\H$, and how they act adds more information.  For it to make
sense to define the action of operators, we need a ``quantum Hilbert space'' $\H$ for them to act on, and not just a function
$Z(q)$.  I explain later how to define natural operators, associated to link cobordisms, that act on $\H$.

\def\I{{\mathcal I}}
The ability to do all this rests on the following facts about the KW equations.  (These facts were also discovered by A. Haydys
\cite{Haydys}.)   I will just state these facts as facts -- which one can verify by a short calculation -- without describing the
quantum field theory construction that motivated me to look for them \cite{WittenK}.    We consider  the KW equations on
a four-manifold $M=W\times \I$, where $W$ is a three-manifold with local coordinates $x^i,\,i=1,2,3$, and $\I$ is a one-manifold
parametrized by $y$.  (In our application, $\I$ will be $\R_+$.   In the first lecture, $\R_+$ was introduced as the direction of
a gradient flow, and parametrized by ``time,'' but  now we interpret $\I$ as a ``space'' direction; we are about to introduce 
 a new ``time''  direction.)   We write $\phi=\sum_{i=1}^3
\phi_i \d x^i+\phi_y\d y$.   Now we replace $M$ by a five-manifold $X=\R\times M$, where  $\R$ is parametrized by a new ``time''
coordinate $t$.  We convert the four-dimensional KW equations on $M$  into five-dimensional equations on $X$ by simply replacing
$\phi_y$, wherever it appears, by a covariant derivative in the new time direction:
\begin{equation}\label{wexo}\phi_y\to \frac{D}{D t}=\frac{\partial}{\partial t}+[A_t,\,\cdot\,]. \end{equation}

If we make this substitution in a random differential equation containing $\phi_y$, we will not get a differential equation but
a differential operator.  In the case of the KW equations, $\phi_y$ appears only inside commutators $[\phi_i,\phi_y]$ and
covariant derivatives $D_\mu\phi_y$, and the substitutions proceed by
\begin{equation}\label{mefo} [\phi_i,\phi_y]\to -D_y\phi_i,~~ D_\mu\phi_y\to [D_\mu, D_t]=F_{\mu t}. \end{equation}
This is enough to show that the substition does give a differential equation.  Generically, the differential equation obtained
this way would not be well-posed, where here well-posed means ``elliptic.''  Essential to make our story work is that
the five-dimensional equation obtained in the case of the KW equations from the substitution $\phi_y\to D/Dt$ actually is elliptic.
 This is  not hard to verify if one suspects it.

The five-dimensional equation has a four-dimensional symmetry that is not obvious from what we have said so far.  We started
on $M=W\times\I$, with $W$ a three-manifold, and then via $\phi_y\to D/Dt$, we replaced $M$ with $X=\R\times M
=\R\times W\times \I$.  It turns out that here $\R\times W$ can be replaced by any oriented four-manifold $Z$, and our
equation can be naturally defined on\footnote{More generally \cite{Haydys}, the equation can be defined on any oriented five-manifold
$X$ endowed with an everywhere nonzero vector field (which for $X=Z\times \I$ we take to be $\partial/\partial y$).}
 $X=Z\times \I$.  At a certain point, we will make use of this four-dimensional symmetry.
 
 Another essential fact is that the five-dimensional equation that we get this way can be formulated as a gradient flow equation
 \begin{equation}\label{melz}\frac{\d \Phi}{\d t}=-\frac{\delta\Gamma}{\delta \Phi}, \end{equation}
 for a certain functional $\Gamma(\Phi)$ (here all fields $A,\phi$ are schematically combined into $\Phi$).  This means that we 
 are in the situation explored by Floer when he defined
 Floer cohomology in the 1980's: modulo analytic subtleties, we can define an infinite-dimensional version of Morse theory,
 with $\Gamma$ as a middle-dimensional Morse function. 
 
 In Morse homology (that is, in Morse theory formulated by counting of gradient flow lines \cite{WittenMorse,Hutchings}),
 one defines a vector space $\V$ with a basis vector $|i\rangle$ for every critical point of $\Gamma$. $\V$ is $\Z$-graded
 by a ``fermion number'' operator $F$ that assigns to $|i\rangle$ the Morse index $f_i$ of the  critical point $i\in\S$.  
 If as in our case, the Morse function is defined on a space that has connected components labeled by another quantity $P$
 (in our case, $P$ is the instanton number operator), then
 $\V$ is also graded by the value of $P$ for a given critical point.    On $\V$, one defines a ``differential'' $Q:\V\to\V$ by 
 \begin{equation}\label{zelb}Q|i\rangle =\sum_{j\in\S|f_j=f_i+1}n_{ij}|j\rangle, \end{equation}
 where for each pair of critical points $i,j$ with $f_j=f_i+1$, we define $n_{ij}$ as the ``number''  of solutions of the
 gradient flow equation 
 \begin{equation}\label{zonel}\frac{\d\Phi}{\d t}=-\frac{\partial\Gamma}{\partial \Phi},~~-\infty<t<\infty \end{equation}
 that start at $i$ in the far past and end at $j$ in the far future.   In the counting, one factors out by the time translation symmetry,
 and one includes a sign given by the sign of the fermion determinant, and is, the sign of the determinant of the linear operator
 obtained by linearizing the gradient flow equation.  On a finite-dimensional compact manifold $B$, the cohomology of $Q$
 is simply the cohomology of $B$ with integer coefficients.  Floer's basic idea (which can be interpreted physically as
 a generalization of the procedure just described to quantum field theories in higher dimension) is that the cohomology
 of $Q$ makes sense in an infinite-dimensional setting provided the flow equation is elliptic and certain compactness properties
 hold.  (In our present context, the flow equation is certainly elliptic but the necessary compactness properties have not yet been
 proved.)
 
\begin{figure}
 \begin{center}
   \includegraphics[width=2in]{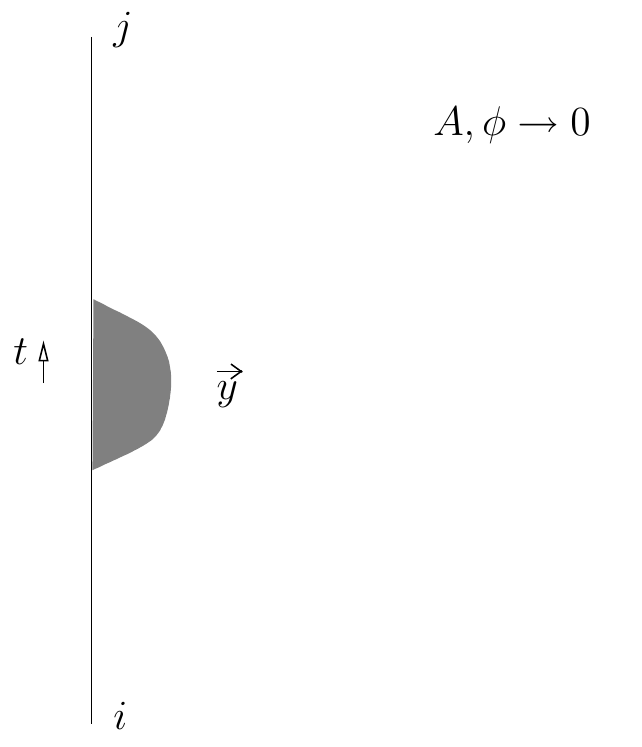}
 \end{center}
\caption{\small   Matrix elements of the differential $Q$ are computed by counting solutions of the five-dimensional
equations interpolating between two given four-dimensional solutions $i$ in the past and $j$ in the future.
The solutions vanish for $y\to\infty$, and at $y=0$ obey boundary conditions that will be described in the text. \label{roughready}}
\end{figure}

 When we follow this recipe  in the present context, the time-independent 
 solutions in five dimensions are just the solutions of the KW equations in
 four dimensions (with $A_t$ reinterpreted as $\phi_y$), since when we ask for a solution to be time-independent, we undo what
 we did to go from four to five dimensions.  So the space $\V$ on which the differential of Morse theory acts is the same space 
 we introduced before in writing the quantum knot invariant $Z(q)$ as a trace.  Moreover, in our application, the procedure
 described in the last paragraph means that the matrix elements $n_{ij}$ of the differential should be computed by
 ``counting'' the five-dimensional solutions that interpolate from a KW solution $i$ in the past to a KW solution $j$ in the future
 (fig. \ref{roughready}).

A conspicuous gap here is that I have not yet described the boundary condition that should be used at the finite end of $X$,
in other words at $y=0$.  Before doing so, I want to describe something interesting that happens in Khovanov homology for
certain gauge groups.  In the study of Khovanov homology for $G^\vee=\SU(2)$, it has been found \cite{ORS} that there are two variants
of the theory, called ``even'' and ``odd'' Khovanov homology.  They are defined using a complex $\V$ that additively is the same
in the two cases, but on this complex one defines two different differentials, say $Q_+$ for the even theory and $Q_-$ for the
odd theory.  They are both defined over $\Z$ and they are congruent mod 2, so their cohomologies, which are called even
and odd Khovanov homology, are isomorphic if one reduces mod 2.  Why would this happen in our framework and for what
groups should we expect it to happen?  I claim that we should use the exceptional isomorphism $\SU(2)\cong\Spin(3)$ and that
in general the bifurcation into even and odd Khovanov homology will occur precisely for $G^\vee=\Spin(2n+1)$, $n=1,2,3,\dots$.

\def\L{{\mathcal L}}
In general, the cohomology of a manifold $B$  can be twisted by a flat complex line 
bundle $\L$.  Instead of the ordinary cohomology $H^i(B,\Z)$, we can consider the cohomology with values in $\L$,
$H^i(B,\L)$.  There is a Morse theory recipe to compute this, by slightly modifying the procedure described above.  The possible
$\L$'s are classified by $\mathrm{Hom}(\pi_1(B),\C^*)$.  In the
present context, $B$ is a function space, consisting of pairs $A,\phi$ on $M=W\times \R_+$ (which represent initial data for
solutions on $X=\R\times M$ where $\R$ is parametrized by ``time,'' and which obey certain boundary conditions).   We only care about the pairs $A,\phi$ up to $G$-valued
gauge transformations (which because of the boundary conditions are trivial on the finite and infinite boundaries of $M$). For
$W=\R^3$, this means that $\pi_1(B)=\pi_4(G)$.  For the simple Lie groups, we have 
\begin{equation}\notag\pi_4(G) = \begin{cases} \Z_2 & G=\Sp(2n) ~{\mathrm {or}} ~\Sp(2n)/\Z_2,~~n\geq 1 \cr
              0 & {\mathrm {otherwise}} .\end{cases}\end{equation} 
So Khovanov homology is unique unless $G^\vee=\SO(2n+1)$, $G=\Sp(2n)$ (or $G^\vee=\Spin(2n+1)$, $G=\Sp(2n)/\Z_2$),
for some $n\geq 1$, in which case there are two versions of Khovanov homology.   Concretely, an $\Sp(2n)$ bundle on a five-dimensional spin manifold
$Y$ (with a trivialization at infinity along $Y$) has a $\Z_2$-valued invariant $\zeta$ derived from $\pi_4(\Sp(2n))=\Z_2$. ($\zeta$ is defined as the mod 2 index
of the Dirac operator valued in the fundamental representation of $\Sp(2n)$.)  When we
define the differential by counting five-dimensional solutions, we have the option to modify the differential by including a factor
of $(-1)^\zeta$.  This gives a second differential $Q'$ that still obeys $(Q')^2=0$, and is congruent mod 2 to the differential $Q$
that is obtained without the factor of $(-1)^\zeta$.  The two theories associated to $Q$ and $Q'$ are the candidates for the
two versions of Khovanov homology.\footnote{Odd Khovanov homology of $SU(2)$ is known \cite{Lauda} to be related to the supergroup
$\mathrm{OSp}(1|2)$.  This connection will be explored elsewhere from a physical point of view \cite{MiW}.}

Next I come to an explanation of the boundary conditions that we
impose on the four- or five-dimensional equations. These boundary
conditions are crucial, since for instance it is only via the boundary
conditions that knots enter. First I will describe the boundary
condition in the absence of knots. It is essentially enough to
describe the boundary condition in four dimensions rather than five
(once one understands it, the lift to five dimensions is fairly
obvious), and as the boundary condition is local, we assume
initially that the boundary of the four-manifold is just $\R^3$. So we
work on $M = \R^3\times \R_+$. (This special case is anyway the right case
for the Jones polynomial, which concerns knots in $\R^3$ or
equivalently $S^3$.)

\def\t{{\frak t}}
As a preliminary to describing the boundary condition, I need to describe an important
equation in gauge theory, which is Nahm's equation.  This is a system of ordinary differential
equations for a triple $X_1,X_2,X_3$ of elements of $\frak g$, the Lie algebra of a compact Lie group $G$.
The equations read
\begin{equation}\label{zolf}\frac{\d X_1}{\d y}+[X_2,X_3]=0, \end{equation}
and cyclic permutations thereof.  On an open half-line $y>0$, Nahm's equations have the special solution
\begin{equation}\label{olf}X_i=\frac{\t_i}{y}, \end{equation}
where $\t_i$ are elements of $\frak g$ that obey the $\frak{su}(2)$ commutation relations $[\t_1,\t_2]=\t_3$ and
cyclic permutations.  Thus the $\t_i$ are the images of a standard basis of $\frak{su}(2)$ under a homomorphism
$\varrho:\frak{su}(2)\to\frak g$.   This  singular solution of Nahm's equations has been important in numerous applications,
in work by Nahm, Kronheimer, Atiyah and Bielawski, and others \cite{Nahm,K,Ktwo,AB}.  
We will use it to define an elliptic boundary condition
on the KW equations and their five-dimensional cousins.  

In fact, Nahm's equations can be embedded in the KW equations (\ref{melf}) on $\R^3\times \R_+$.  We look for a solution that
{\it (i)} is invariant under translations of $\R^3$, {\it (ii)} has the property that $A=0$, and {\it (iii)} has the property that
$\phi=\sum_{i=1}^3\phi_i \,\d x^i+0\cdot \d y$.  For solutions satisfying these conditions, the KW equations reduce to Nahm's equations
\begin{equation}\label{worff}\frac{\d \phi_1}{\d y}+[\phi_2,\phi_3]=0, \end{equation}
and cyclic permutations.  So the ``Nahm pole'' gives a special solution
\begin{equation}\label{orf} \phi_i = \frac{\t_i}{y}.\end{equation}
We define a boundary condition by saying that we only allow solutions that are asymptotic to this one for $y\to 0$, modulo
less singular terms.    See \cite{MW} for a proof that, for any $\varrho$,
 this boundary condition is elliptic and has regularity properties similar to those
of more standard elliptic boundary conditions such as Dirichlet or Neumann.  (Actually, in that paper, a somewhat wider class of
 Nahm pole boundary conditions is analyzed.)  For applications to the Jones polynomial and Khovanov
homology, we take $\varrho:\frak{su}(2)\to \frak g$ to be a principal embedding in the sense of Kostant (for example, 
for $G=SU(n)$, this
means that the $n$-dimensional representation of $G$ transforms as an irreducible representation of $\frak{su}(2)$).  

This is the boundary condition that we want at $y=0$, in the absence of knots.  For the most simple application to the quantum
knot invariants and Khovanov homology, we require that $A,\phi\to 0$ for $y\to\infty$.  With these conditions at $y=0,\infty$,
and suitable conditions for $x^i\to\infty$, it is possible to prove \cite{MW} that the Nahm pole solution is the only solution
on $\R^3\times \R_+$.  This corresponds to the statement that the Khovanov homology of the empty link is of rank 1.

   \begin{figure}
 \begin{center}
   \includegraphics[width=3in]{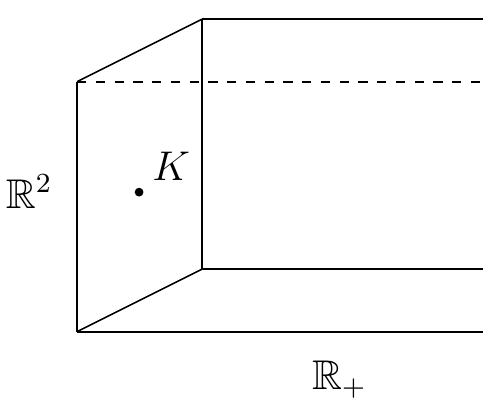}
 \end{center}
\caption{\small   The boundary condition in the presence of a knot $K$ is determined by a model solution of
a reduced equation on $\R^2\times \R_+$.  Near a generic boundary point, the model solution has the Nahm pole
singularity, but it has a more complicated singular behavior near the boundary point corresponding to $K$.  This
model solution depends on the choice of a representation $R^\vee$ of the dual group $G^\vee$.   \label{reduced}}
\end{figure}

Now I should explain how the boundary condition is modified along a knot $K$.  This will be done by requiring a more subtle
singularity along the knot.  The local model is that the boundary of $M$ is $\R^3$ and the knot $K$ is a copy of $\R\subset \R^3$.
The boundary condition is defined by giving a model solution on $\R^3\times \R_+$ that away from $K$ has the now familiar
Nahm pole singularity at $y=0$, but has a more complicated singular behavior along $K$.  The model solution is invariant under
translations along $K$, so it can be obtained by solving some reduced equations on $\R^2\times \R_+$.  In the reduced
picture, $K$ corresponds to a point in $\R^2\times \{y=0\}$ (fig. \ref{reduced}).  The model solution depends on the choice
of a representation $R^\vee$ of the dual group $G^\vee$.  Near a boundary point disjoint from $K$, the model solution
has the usual Nahm pole singularity, but near $K$ it has a more complicated singularity.  The relevant model solutions
can be found in closed form.  This has been done in \cite{WittenK}, section 3.6, for $G^\vee=\SU(2)$, and in
\cite{mikhaylov} for any $G^\vee$.   I will not describe these reduced solutions here, except to say that the singularity
along $K$ is a more complicated cousin of the singularity used in \cite{KW} to describe the geometric Hecke transformations
of the geometric Langlands correspondence.

    \begin{figure}
 \begin{center}
   \includegraphics[width=2in]{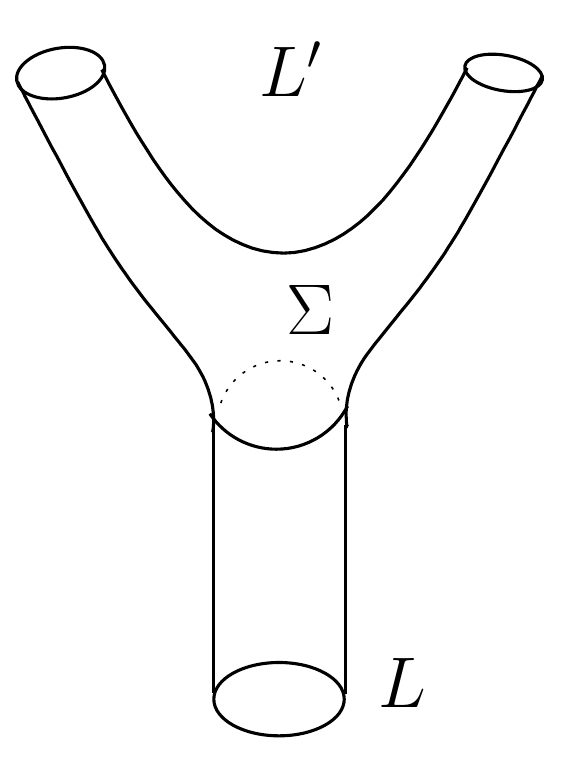}
 \end{center}
\caption{\small   A two-dimensional surface $\Sigma$ that represents a ``link cobordism'' from a link $L$ in the past to
another link $L'$ in the future.  In the example shown, for simplicity, $L$ is an unknot and $L'$ consists of two unlinked unknots.
$\Sigma$ is embedded in the boundary of $X=\R^4\times\R_+$.  \label{Surface}}
\end{figure}

   \begin{figure}
 \begin{center}
   \includegraphics[width=3.5in]{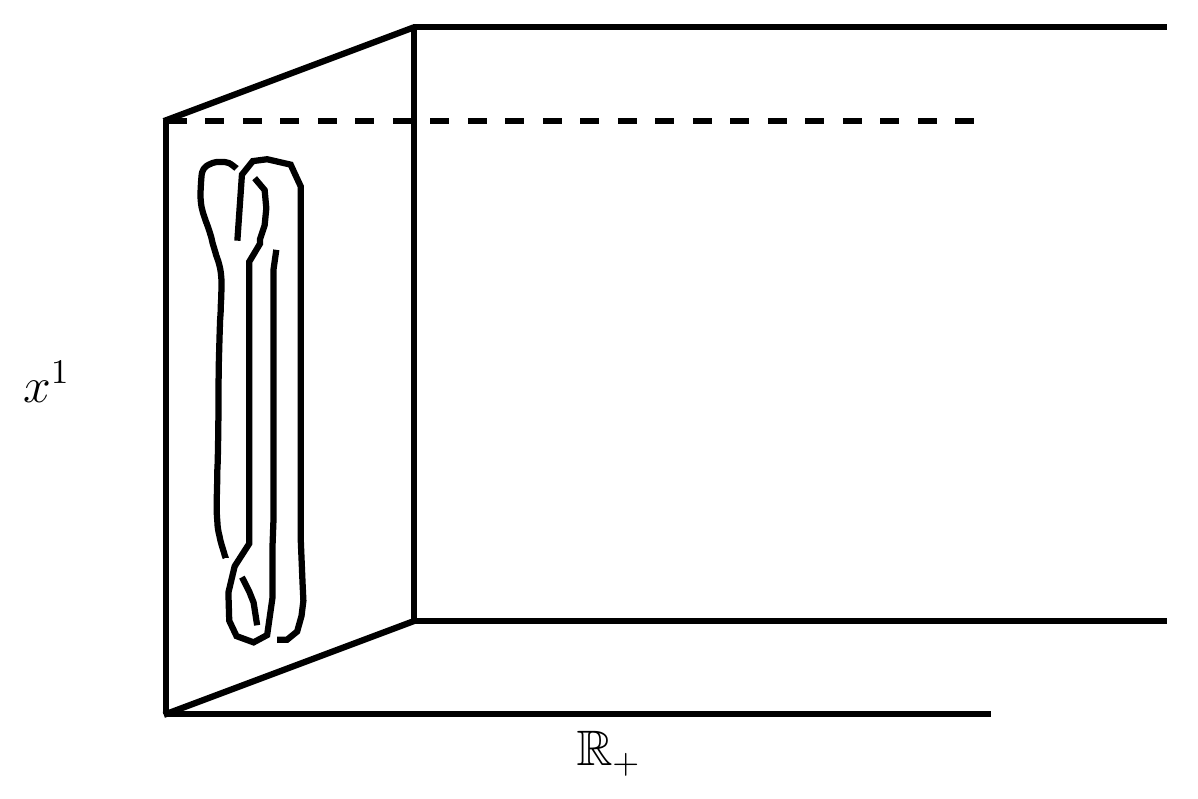}
 \end{center}
\caption{\small Stretching a knot in one dimension, to reduce to a description in one dimension less.  One of the
directions in $\R^3$ -- here labeled as $x^1$ -- plays the role of ``time.''  After stretching of the knot,
one hopes that a solution of the equations becomes almost everywhere nearly independent of $x^1$. If so,
a knowledge of the solutions that are actually independent of $x^1$ can be a starting point for understanding
four-dimensional solutions.  This type of analysis will fail at the critical points of the function $x^1$ along a knot or link;
a correction has to be made at those points.}
 \label{stretch}
\end{figure} 
 The model solution has a singularity that, along the boundary, is of codimension 2.  When we go to five dimensions, the
singularity remains of codimension 2, so now, as the boundary is a four-manifold, the singularity is supported on a two-dimensional
surface $\Sigma$, not a knot (fig. \ref{Surface}).  A boundary condition modified on a 2-surface in the boundary is what we need to
define the ``morphism'' of Khovanov homology associated to a ``link cobordism.''  In other words, given a 2-surface $\Sigma$
that interpolates between one link $L$ in the past and another link $L'$ in the future, as sketched in the figure,
counting solutions with boundary conditions
modified along $\Sigma$ gives the matrix element for a time-dependent transition from a physical state (a cohomology class of
$Q$) in the past in the presence of  $L$ to a physical state in the future in the presence of $L'$.  The morphisms that are defined this
way have the formal properties that are expected in Khovanov homology.

There is another reason that it is important to describe the
reduced equations in three dimensions. To compute the Jones
polynomial, we need to count certain solutions in four dimensions;
knowledge of these solutions is also the first step in constructing
the candidate for Khovanov homology. How are we supposed to
describe four-dimensional solutions? A standard strategy, often
used in Floer theory and its cousins, involves ÒstretchingÓ the knot
in one direction, in the hope of reducing to a piecewise description
by solutions in one dimension less (fig. \ref{stretch}).  To get anywhere with such an analysis,
we need to be able to solve the three-dimensional reduced equations.
  
 \begin{figure}[ht]
 \begin{center}
   \includegraphics[width=2.2in]{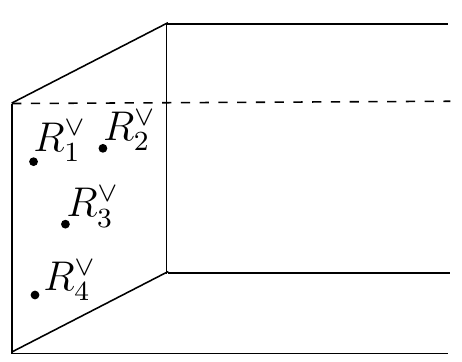}
 \end{center}
\caption{\small $\R^2\times\R_+$ with $n$ boundary points labeled by representations $R_1^\vee,\dots, R_n^\vee$ of the dual group, in this case with $n=4$.}
 \label{labeling}
\end{figure} 
  
 Another way to make the point is as follows. Most mathematical
definitions of Khovanov homology proceed, directly or implicitly, by
defining a category of objects associated to a two-sphere (or in
some versions, a copy of $\C=\R^2$) with marked points that are
suitably labeled.
In the present approach, this category should be a category associated to
 solutions of the reduced
three-dimensional equations on $\R^2\times \R_+$, with finitely many points
 in $\R^2\times \{y=0\}$ labeled by representations of the dual group $G^\vee$
 (fig. \ref{labeling}).  By analyzing this situation, Gaiotto and I \cite{GWagain} were able
 to get a fairly clear framework for understanding the relevant category. (We found that to make this program
 work nicely, we had to perturb to a slightly more generic version of the problem than I have described.)   For $G^\vee=\SU(2)$,
 the category is expected to be a Fukaya-Seidel category (an $A$-model category with a superpotential) where the target
 space is a moduli space of monopoles on $\R^3$, and the  superpotential  encodes the positions of
 the knots.  This category has not yet been analyzed in any detail,   but for the more modest problem of understanding the Jones
 polynomial (rather than Khovanov homology), we were able to get a reasonably satisfactory picture.

 In doing this, we used
 a simplification that can be achieved by modifying the condition on how a solution of the KW equations should
 behave for $y\to \infty$.  We kept the condition that $A\to 0$ at infinity, but required $\phi$ to approach
 $\sum_{i=1}^3 c_i\d x^i$, where the $c_i$ are a prescribed triple of commuting elements of $\frak g$. 
 (It is known \cite{Ktwo} that this is a convenient asymptotic condition in the study of Nahm's equations.)
 The counting of solutions of an elliptic differential equation is invariant under continuous deformations, as long
 as certain conditions are satisfied. So it is expected that the counting of solutions of the KW equations is independent of
 $\vec c=(c_1,c_2,c_3)$.  By exploiting a suitable choice of $\vec c$, we were able to relate the counting of solutions of the KW
 equations to the vertex model for the Jones polynomial, which was the starting point in Lecture One.  Notice that $\vec c$
 defines a direction in $\R^3$.  This direction determines the knot projection that is used in the vertex model.

Recent work reconsidering the Fukaya-Seidel category from a physical point of view \cite{GMW} may be helpful
in understanding better the  categories that arise in the present context when a knot is stretched in one
dimension.

Research supported in part by NSF Grant PHY-1314311.

\bibliographystyle{unsrt}

\begin{thebibliography}{99}


\bibitem{Jones}
V. F. R. Jones, ``A Polynomial Invariant For Links Via Von Neumann Algebras,'' Bull.
AMS {\bf{12}} (1985) 103.

\bibitem{Jonestwo}
V. F. R. Jones, ``On Knot Invariants Related To Some Statistical 
Mechanical Models,'' Pacific J. Math. {\bf 137} (1989) 311-334.

\bibitem{Kauffman}
L. H. Kauffman, {\it Knots And Physics} (World Scientific, 1991).

\bibitem{TK}
A. Tsuchiya and Y. Kanie, ``Vertex Operators In Conformal Field Theory On
$\Bbb P^1$ And Monodromy Representations of Braid Group,'' Adv. Stud. Pure
Math. {\bf 16} 297-372.

\bibitem{WittenJones}
E. Witten, ``Quantum Field Theory And The Jones Polynomial,'' Commun. Math. Phys. {\bf 121}
(1989) 351-399.

\bibitem{FK}
J. Frohlich and C. King, ``The Chern-Simons Theory And Knot Polynomials,'' Commun. Math. Phys. {\bf 126} (1989) 167-199,

\bibitem{Kash} R. M. Kashaev, ``The Hyperbolic Volume Of Knots
From The Quantum Dilogarithm,'' Lett. Math. Phys. {\bf 39} (1997)
269-275.

\bibitem{MM} H. Murakami and J. Murakami, ``The Colored
Jones Polynomial And The Simplicial Volume Of A Knot,'' Acta Math.
{\bf 186} (2001) 85-104.

\bibitem{MMOTY}
 H. Murakami, J. Murakami, M. Okamoto, T. Takata, and Y. Yokota,
 ``Kashaev's Conjecture And The Chern-Simons Invariants Of Knots
 And Links,'' Experiment. Math. {\bf 11} (2002) 427-435.

\bibitem{KT}
R. M. Kashaev and O. Tirkkonen, ``A Proof Of The Volume Conjecture
For Torus Knots,'' J. Math. Sci. (NY) {\bf 115} (2003) 2033-2036,
math.GT/9912210.

\bibitem{Gu}
S. Gukov, ``Three-Dimensional Quantum Gravity,
Chern-Simons Theory, And The $A$-Polynomial,'' Commun. Math. Phys.
{\bf 255} (2005) 577-627.

\bibitem{Mm}
H. Murakami, ``Asymptotic Behaviors Of The Colored Jones
Polynomials Of A Torus Knot,'' Internat. J. Math. {\bf 15} (2004)
547-555.

\bibitem{Analytic}
E. Witten, ``Analytic Continuation Of The Jones Polynomial,'' arXiv:1001.2933.

\bibitem{FG}
D. Freed and R. Gompf, ``Computer Calculations Of Witten's
3-Manifold Invariant,'' Commun. Math. Phys. {\bf 141} (1991)
79-117.

\bibitem{LJe}
L. Jeffrey, ``Chern-Simons-Witten Invariants Of Lens Spaces And
Torus Bundles and the Semi-Classical Approximation,'' Commun.
Math. Phys. {\bf 147} (1992) 563-604.

\bibitem{Wittengauge}
E. Witten, ``Khovanov Homology and Gauge Theory,''  in  R. Kirby, V. Krushkal, and Z. Wang, eds.,
{\it Proceedings Of The FreedmanFest} (Mathematical Sciences Publishers, 2012) 291-308, arXiv:1108.3103.

 \bibitem{KW}
 A. Kapustin and E. Witten,  ``Electric-Magnetic Duality And The Geometric Langlands Program,'' Commun. Numb. Th. Phys. {\bf 1} (2007) 
1-236,  hep-th/0604151.

\bibitem{Taubes}
C. H. Taubes,  ``Compactness Theorems For $SL(2;\C)$ Generalizations Of The Anti-Self Dual Equations, Part I,''
arXiv:1307.6447.

\bibitem{Taubestwo}
C. H. Taubes,  ``Compactness Theorems For $SL(2;\C)$ Generalizations Of The Anti-Self Dual Equations, Part II,''
arXiv:1307.6451.


\bibitem{GU}
M. Gagliardo and K. Uhlenbeck, ``The Geometry Of The Kapustin-Witten Equations,'' J. Fixed Point Theory Appl. {\bf 11} (2012) 185-198.


\bibitem{WittenK}
E. Witten, ``Fivebranes And Knots,'' Quantum Topology {\bf 3} (2012) 1-137, arXiv:1101.3216.

\bibitem{GW}
D. Gaiotto and E. Witten, ``Supersymmetric Boundary Conditions In $\N=4$ Super Yang-Mills Theory,'' J. Stat. Phys. {\bf 135} (2009) 789-855,
 arXiv:0804.2902.


\bibitem{AtiyahWeyl}
M. F. Atiyah, ``New Invariants Of Three And Four DImensional Manifolds,'' in {\it The Mathematical Heritage Of Herman Weyl},
Proc. Symp. Pure Math. {\bf 48} (1988) 285-300.
 

 
 \bibitem{OV}
H. Ooguri and C. Vafa, ``Knot Invariants And Topological Strings,''  Nucl. Phys.
{\bf{B577}}
(2000) 419, hep-th/9912123.
 
 \bibitem{GV}
R. Gopakumar and C. Vafa, ``On The Gauge Theory/Geometry Correspondence,'' Adv. Theor.
Math. Phys. {\bf 3} (1999) 1415-1443, hep-th/9811131.



\bibitem{GV2}
R. Gopakumar and C. Vafa, ``$M$-Theory And Topological Strings, I, II,''
hep-th/9809187,
hep-th/9812127.
 
 \bibitem{Khovanov}
M. Khovanov, ``A Categorification Of The Jones Polynomial,'' Duke. Math. J. {\bf
101} (2000)
359-426.
 
 \bibitem{GSV}
S. Gukov, A. S. Schwarz, and C. Vafa,
``Khovanov-Rozansky Homology And Topological Strings,''
Lett. Math. Phys. {\bf 74} (2005) 53-74, hep-th/0412243.
 
 \bibitem{CK}
S. Cautis and J. Kamnitzer, ``Knot Homology Via Derived Categories Of Coherent
Sheaves I,
$\frak{sl}(2)$ Case,'' arXiv:math/0701194.

\bibitem{SS}
P. Seidel and I. Smith, ``A Link Invariant From The Symplectic Geometry Of Nilpotent
Slices,''
arXiv:math/0405089.

\bibitem{Kam}
J. Kamnitzer, ``The Beilinson-Drinfeld Grassmannian And Symplectic Knot Homology,''
arXiv:0811.1730.

\bibitem{KM}
P. B. Kronheimer and T. S. Mrowka,  ``Knot Homology Groups From Instantons,''
arXiv:0806.1053.

\bibitem{KM2}
P. B. Kronheimer and T. S. Mrowka, ``Khovanov Homology Is An Unknot-Detector,''
arXiv:1005.4346.   

\bibitem{Haydys}
A. Haydys, ``Fukaya-Seidel Category And Gauge Theory,''
arXiv:1010.2353.

\bibitem{WittenMorse}
  E.~Witten,
  ``Supersymmetry and Morse Theory,''
  J.\ Diff.\ Geom.\  {\bf 17} 661 (1982).

\bibitem{Hutchings}
M. Hutchings, ``Lecture Notes On Morse homology (with an Eye Towards Floer Homology and Pseudoholomorphic Curves),''
available at math.berkeley.edu/hutching.

\bibitem{ORS}
P. Ozsvath, J. Rasmussen, and Z. Szabo, ``Odd Khovanov Homology,'' arXiv:0710.4300.

\bibitem{Lauda}
A. P. Ellis and A. D. Lauda, ``An Odd Categorification
Of Quantum $sl(2)$,'' arXiv:1307.7816.


\bibitem{MiW}
V. Mikhalylov and E. Witten, ``Branes and Supergroups,'' to appear.

\bibitem{Nahm}
W. Nahm, ``A Simple Formalism For The BPS Monopole,'' Phys. Lett. {\bf B90} (1980) 413.

\bibitem{K}
P. Kronheimer, ``Instantons And The Geometry Of The Nilpotent Variety,'' J. Diff. Geom.
{\bf 32} (1990) 473-90.

\bibitem{Ktwo}
P. Kronheimer, ``A Hyper-Kahlerian Structure On Coadjoint Orbits Of A Semisimple Complex Lie Group,''
J. London Math. Soc. {\bf 42} 193-206 (1990).

\bibitem{AB}
M. F. Atiyah and R. Bielawski, ``Nahm's Equations, Configuration Spaces, and Flag Manifolds,'' 
 arXiv:math/0110112.

\bibitem{MW}
R. Mazzeo and E. Witten, ``The Nahm Pole Boundary Condition,'' arXiv:1311.3167.

\bibitem{mikhaylov}
V. Mikhaylov, ``On The Solutions Of Generalized Bogomolny Equations,''  arXiv:1202.4848.

\bibitem{GWagain}
D. Gaiotto and E. Witten, ``Knot Invariants From Four-Dimensional Gauge Theory,'' arXiv:1106.4789.

\bibitem{GMW}
D. Gaiotto, G. W. Moore, and E. Witten, to appear.
\end{thebibliography}

\end{document}